\documentclass[12pt]{article}
\usepackage{amsmath,amssymb,amsthm, amscd}

\begin{document}

\begin{center}
{\bf A characterization of nilpotent orbit closures among symplectic singularities}
\end{center}
\vspace{0.4cm}

\begin{center}
{\bf Yoshinori Namikawa}
\end{center}
\vspace{0.2cm}

Symplectic singularities have been playing important roles both in algebraic geometry and 
geometric representation theory ever since Beauville introduced their notion in [Be]. Most 
examples of symplectic singularities admit natural $\mathbf{C}^*$-actions with 
only positive weights. Kaledin [Ka] conjectured that any symplectic singularity admits such a $\mathbf{C}^*$-action. 

If a symplectic singularity has a $\mathbf{C}^*$-action with positive weights, it can be globalized to an affine variety with a $\mathbf{C}^*$-action. Such an affine variety is called a {\em conical symplectic variety}. More precisely, an affine normal variety $X = \mathrm{Spec}\; R$ is a conical symplectic variety if 

(i) $R$ is positively graded: $R = \oplus_{i \ge 0}R_i$ with $R_0 = \mathbf{C}$; 

(ii) the smooth part $X_{reg}$ admits a homogeneous symplectic 2-form $\omega$ and 
it extends to a holomorphic 2-form on a resolution $\tilde{X}$ of $X$. 

Denote the $\mathbf{C}^*$-action by $t: X \to X$ $(t \in \mathbf{C}^*)$. By the assumption we have $t^*\omega = t^{l}\omega$ for some integer $l$. This integer $l$ is called the weight of $\omega$ and is denoted by $wt(\omega)$. By the extension property (ii) we have $wt(\omega) > 0$ (cf. [Na 3], Lemma (2.2)).  

Let $\{x_0, ..., x_n\}$ be a set of minimal homogeneous generators of the $\mathbf{C}$-algebra $R$ and put $a_i := \mathrm{deg} \; x_i$. We put $N:= \mathrm{max} \{a_0, ..., a_n\}$ and call $N$ the {\em maximal weight} of $X$. It is uniquely determined by a conical symplectic variety $X$. 
By [Na 1], there are only finitely many conical symplectic varieties $(X, \omega)$ of a fixed dimension $2d$ and with a fixed maximal weight $N$, up to an isomorphism. 
In this sense it would be important to classify conical symplectic varieties with maximal weight $1$. 
By the homogeneous generators $\{x_i\}$, we can embed $X$ into an affine space $\mathbf{C}^{n+1}$. 
In [Na 2] we treat the case where $X \subset \mathbf{C}^{n+1}$ is a complete 
intersection of homogeneous polynomials. The main theorem of [Na 2] asserts that $(X, \omega)$ is isomorphic to the nilpotent cone $(N, \omega_{KK})$ of a complex semisimple Lie algebra $\mathfrak{g}$ together with the Kirillov-Kostant 2-form provided that $X$ is singular.  
However, there are a lot of examples of maximal weight $1$ which are not of complete intersection. In fact, a nilpotent orbit $O$ of a complex semisimple Lie algebra $\mathfrak{g}$ admits the  Kirillov-Kostant form $\omega_{KK}$ and if its closure $\bar{O}$ is normal, then  $(\bar{O}, \omega_{KK})$ is a conical symplectic variety with maximal weight $1$ by  Panyushev [Pa] and Hinich [Hi]. 

A main purpose of this article is to prove that they actually exhaust all conical symplectic varieties with maximal weight $1$.     
\vspace{0.2cm}

{\bf Theorem}. {\em Let $(X, \omega)$ be a conical symplectic variety with maximal weight  $1$. Then $(X, \omega)$ is isomorphic to one of the following:} 

(i) {\em $(\mathbf{C}^{2d}, \omega_{st})$ with $\omega_{st} = \Sigma_{1 \le i \le d} dz_i \wedge dz_{i + d}$}, 

(ii){ \em $(\bar{O}, \omega_{KK})$ where $\bar{O}$ is a normal nilpotent orbit closure of a complex semisimple Lie algebra $\mathfrak{g}$ and $\omega_{KK}$ is the Kirillov-Kostant form.}
\vspace{0.2cm}

There is a non-normal nilpotent orbit closure in a complex semisimple Lie algebra. The  normalization $\tilde{O}$ of such an orbit closure $\bar{O}$ is also a conical symplectic variety.\footnote{By [K-P, Proposition 7.4] $\bar{O}$ is always resolved by a vector bundle $Y$ over $G/P$ with a parabolic subgroup $P$ of the adjoint group $G$ of $\mathfrak{g}$. Denote this resolution by $\pi : Y \to \bar{O}$. The map $\pi$ factorizes as $Y \to \tilde{O} \to \bar{O}$. The fiber $\pi^{-1}(0)$ coincides with the zero section of $Y$, which is isomorphic to $G/P$. As $G/P$ is connected, the fibre $\mu^{-1}(0)$ of the normalization map $\mu : \tilde{O} \to \bar{O}$ consists of just one point, say $x \in \tilde{O}$. The $\mathbf{C}^*$-action on $\bar{O}$ extends to a $\mathbf{C}^*$-action on $\tilde{O}$ with a unique fixed point $x$. 
It is easily checked that this $\mathbf{C}^*$-action has only positive weights and $\tilde{O}$ 
becomes a conical symplectic variety. }
But the maximal weight of $\tilde{O}$ is usually larger than $1$.\footnote{It may happen that $\tilde{O}$ coincides with a normal nilpotent orbit closure of a different complex 
semisimple Lie algebra. In such a case the maximal weight is $1$.}  

We first notice that $\omega$ determines a Poisson structure on $X_{reg}$ in a usual way. 
By the normality of $X$, it uniquely extends to a Poisson structure $\{\cdot , \; \cdot \}: 
O_X \times O_X \to O_X$. In particular, $R$ becomes a Poisson $\mathbf{C}$-algebra with 
a Poisson bracket of degree $-wt(\omega)$.  

In the remainder, $X$ is a conical symplectic variety with the maximal weight $N = 1$. 
First of all, we prove in Proposition 1 that $wt(\omega) = 2$ or $wt(\omega) = 1$. In the first case $(X, \omega)$ is isomorphic to an affine space $\mathbf{C}^{2d}$ together with the standard symplectic form $\omega_{st}$.  In the second case the Poisson bracket has degree $-1$ and $R_1$ has a natural Lie algebra structure.  Then it is fairly easy to show that $X$ is a coadjoint orbit closure of a complex Lie algebra $\mathfrak{g}$ (Proposition 3). If $X$ has a crepant resolution, we can prove that $\mathfrak{g}$ is semisimple in the same way as in [Na 2]. But $X$ generally does not have 
such  a resolution and we need a new method to prove the semisimplicity. This is nothing but  Proposition 4. 
\vspace{0.2cm}

{\bf Proposition 1} {\em Assume that $X$ is a conical symplectic variety with maximal weight $N = 1$. Then $wt(\omega) = 1$ or $wt(\omega) = 2$. If  $wt(\omega ) = 2$, then $(X, \omega)$ is isomorphic to an affine space $(\mathbf{C}^{2d}, \omega_{st})$ with the standard symplectic form..} 

{\bf Remark}. As is remarked in the beginning of [Na 2, \S 2], if $X$ is a smooth conical symplectic variety with maximal weight $1$, then $(X, \omega) \cong (\mathbf{C}^{2d}, \omega_{st})$. 
Hence $X$ is singular exactly when $wt(\omega) = 1$.  

{\em Proof}. Since $N = 1$, the coordinate ring $R$ is generated by $R_1$. 
We put $l := wt(\omega)$. We already know that $l > 0$. 
If $l > 2$, then we have $\{R_1, R_1\} = 0$ and hence $\{R, R\} = 0$, which is absurd. 
We now assume that $l = 2$ and prove that $X$ is an affine space with the standard symplectic form. Then the Poisson bracket induces a skew-symmetric form $R_1 \times R_1 \to R_0 = \mathbf{C}$. If this is a degenerate skew-symmetric form, then we 
can choose a non-zero element $x_1 \in R_1$ such that $\{x_1, \cdot\} = 0$. 
Notice that  
$x_1 = 0$ determines a non-zero effective divisor $D$ on $X_{reg}$. If we choose a general 
point $a \in D$, then the reduced divisor $D_{red}$ is smooth around $a$. Consider an analytic open neighborhood $U \subset X_{reg}$ of $a$. 
Then there is a system of local coordinates $\{z_1, ..., z_{2d}\}$ of $U$ such that  $x_1$ can be  written as $x_1 = z_1^m$ for a suitable $m > 0$. The Poisson structure on $X$ induces a 
non-degenerate Poisson structure $\{\cdot, \cdot \}_U$ on $U$. But, by the choice of 
$x_1$, we have $\{z_1^m, \cdot\}_U = mz_1^{m-1}\{z_1, \cdot \}_U = 0$, which implies that 
$\{z_1, \cdot \}_U = 0$. This contradicts that the Poisson bracket $\{\cdot, \cdot \}_U$ is non-degenerate. 

Therefore the skew-symmetric form is non-degenerate. In this case $X$ is a closed 
Poisson subscheme of an affine space with a non-degenerate Poisson structure induced by 
the standard symplectic form. But such an affine space has no Poisson closed subscheme 
except the affine space itself. Therefore $X = \mathrm{Spec} R$.  Q.E.D.     
\vspace{0.2cm}

The regular part $X_{reg}$ of a conical symplectic variety $X$ is a smooth Poisson variety. 
Let $\Theta_{X_{reg}}$ denote the sheaf of vector fields on $X_{reg}$. 
By using the Poisson bracket we define the Lichnerowicz-Poisson complex 
$$ 0 \to \Theta_{X_{reg}} \stackrel{\delta_1}\to \wedge^2 \Theta_{X_{reg}} \stackrel{\delta_2}\to  
... $$ 
by 
$$ \delta_p f(da_1 \wedge ... \wedge da_{p+1})  
:= \sum_{i = 1}^{p+1} (-1)^{i+1}\{a_i, f(da_1 \wedge 
... \hat{da_i} \wedge ... \wedge da_{p+1}\} $$   
$$ + \sum_{j < k} (-1)^{j+k}f(d\{a_j, a_k\} \wedge da_1 \wedge ... \wedge 
\hat{da_j}  
\wedge ... \wedge \hat{da_k} \wedge ... \wedge da_{p+1}).$$ 
In the Lichnerowicz-Poisson complex, $\wedge^p \Theta_{X_{reg}}$ is placed in degree $p$. 
The Lichnerowicz-Poisson complex  of $X_{reg}$ is closely 
related to the Poisson deformation of $(X, \{\; , \; \})$. For details, see [Na 4]. 

In the remainder we assume that $wt(\omega) = 1$.  
The Poisson bracket then defines a pairing map $R_1 \times R_1 \to R_1$ and $R_1$  becomes a Lie algebra.  We denote by 
$\mathfrak{g}$ this Lie algebra. As all generators have weight 1, we have a surjection 
$\oplus \mathrm{Sym}^{i}(R_1) \to R$. It induces a $\mathbf{C}^*$-equivariant closed embedding $X \to \mathfrak{g}^*$.  

Recall that the adjoint group $G$ of $\mathfrak{g}$ (cf. [Pro, p.86])
is defined as a subgroup of $GL(\mathfrak{g})$ generated by all elements of the form 
$\mathrm{exp}(ad \; v)$ with $v \in \mathfrak{g}$. The adjoint group $G$ is a complex Lie  subgroup of $GL(\mathfrak{g})$, but it is not necessarily a closed algebraic subgroup of 
$GL(\mathfrak{g})$. Moreover, the Lie algebra $Lie(G)$ does not necessarily coincide with 
$\mathfrak{g}$. We have $Lie(G) = \mathfrak{g}$ if and only if the adjoint representation is a 
faithful $\mathfrak{g}$-representation, or equivalently, $\mathfrak{g}$ has trivial center.

\vspace{0.2cm}

{\bf Proposition 2}. {\em Let $\mathrm{Aut}^{{\mathbf C}^*}(X, \omega)$ denote the $\mathbf{C}^*$-equivariant automorphism group preserving $\omega$. Then the identity component of $\mathrm{Aut}^{\mathbf{C}^*}(X, \omega)$ can be identified with the adjoint group $G$ of  $\mathfrak{g}$. Moreover $\mathfrak{g}$ has trivial center.  In particular, $\mathfrak{g}$ is the Lie algebra of the linear algebraic group $\mathrm{Aut}^{\mathbf{C}^*}(X, \omega)$. }   

{\em Proof }. 
Let $(\wedge^{\geq 1}\Theta_{X_{reg}}, \delta)$ be the Lichnerowicz-Poisson complex for the  smooth Poisson variety $X_{reg}$. The algebraic torus $\mathbf{C}^*$ acts on $\Gamma(X_{reg}, \wedge^p\Theta_{X_{reg}})$ and there is an associated grading 
$$ \Gamma(X_{reg}, \wedge^p\Theta_{X_{reg}}) = \oplus_{n \in \mathbf{Z}} \;  
\Gamma(X_{reg}, \wedge^p\Theta_{X_{reg}})(n).$$ 

Since the Poisson bracket of $X$ has degree $-1$, the coboudary map $\delta$ has degree $-1$; thus we have a complex 
$$ \Gamma(X_{reg}, \Theta_{X_{reg}})(0) \stackrel{\delta_1}\to 
\Gamma(X_{reg}, \wedge^2\Theta_{X_{reg}})(-1) \stackrel{\delta_2}\to ... 
$$

The kernel $\mathrm{Ker}(\delta_1)$ of this complex is isomorphic to the tangent space 
of $\mathrm{Aut}^{{\mathbf C}^*}(X, \omega)$ at $[id]$. In fact, an element of 
$\mathrm{Ker}(\delta_1)$ corresponds to a derivation of $O_{X_{reg}}$ (or an infinitesimal automorphism of $X_{reg}$) preserving the Poisson structure, but it uniquely extends to a derivation of $O_X$ preserving the Poisson structure (cf. [Na 4, Proposition 8]).  
   
The Lichnerowicz-Poisson complex $(\wedge^{\geq 1}\Theta_{X_{reg}}, \delta)$ is identified 
with the truncated De Rham complex $(\Omega_{X_{reg}}^{\geq 1}, d)$ by the symplectic 
form $\omega$ (cf. [Na 4, Proposition 9], [Na 3, Section 3]). 
The algebraic torus $\mathbf{C}^*$ acts on 
$\Gamma(X_{reg}, \Omega^p_{X_{reg}})$ and there is an associated grading 
$$ \Gamma(X_{reg}, \Omega^p_{X_{reg}}) = \oplus_{n \in \mathbf{Z}} \; 
\Gamma(X_{reg}, \Omega^p_{X_{reg}})(n).$$
The coboundary map $d$ has degree $0$; thus we have a complex 
$$ \Gamma(X_{reg}, \Omega^1_{X_{reg}})(1) \stackrel{d_1}\to 
\Gamma(X_{reg}, \Omega^2_{X_{reg}})(1) \stackrel{d_2}\to ...  $$ 
Since $\omega$ has weight $1$, this complex is identified with 
the the Lichnerowicz-Poisson complex above. 

There is an injective map $d: \Gamma(X_{reg}, O_{X_{reg}})(1) \to 
\Gamma(X_{reg}, \Omega^1_{X_{reg}})(1)$. 

We shall prove that $\mathrm{Ker}(d_1) = \Gamma(X_{reg}, O_{X_{reg}})(1)$.     
The $\mathbf{C}^*$-action on $X$ defines a vector field 
$\zeta$ on $X_{reg}$.   
For $v \in \Gamma(X_{reg}, \Omega^1_{X_{reg}})(1)$, the Lie derivative $L_{\zeta}v$ of $v$ along $\zeta$
equals $v$. If moreover $v$ is $d$-closed, then one has $v = d(i_{\zeta}v)$ by 
the Cartan relation 
$$ L_{\zeta}v = d(i_{\zeta}v) + i_{\zeta}(dv). $$  
This means that $v \in \Gamma(X_{reg}, O_{X_{reg}})(1)$. 
On the other hand, we have $\Gamma(X_{reg}, O_{X_{reg}})(1) = \Gamma(X, O_X)(1)
= R_1 = \mathfrak{g}$. 

It follows from the identification of $\mathrm{Ker}(\delta_1)$ and $\mathrm{Ker}(d_1)$ that every element of $\mathrm{Ker}(\delta_1)$ is a Hamiltonian vector field $H_f := \{f, \cdot \}$ for some $f \in R_1$. In particular, for $g \in R_1$, we have $H_f(g) = [f, g]$. Since $H_f \ne 0$ for a non-zero $f$, the map $ad: \mathfrak{g} \to \mathrm{End}(\mathfrak{g})$ is an injection. 
Notice that an element of $\mathrm{Aut}^{\mathbf{C}^*}(X, \omega)$ determines an automorphism of a graded $\mathbf{C}$-algebra $R$. In particular, it induces a $\mathbf{C}$-linear automorphism of $R_1 = \mathfrak{g}$. Since $R$ is generated by $R_1$, this linear automorphism completely determines an automorphism of $R$. Hence,     
both $G$ and $\mathrm{Aut}^{\mathbf{C}^*}(X, \omega)$ are subgroups of 
$GL(\mathfrak{g})$. The tangent spaces of both subgroups at $[id]$ coincide with  $\mathfrak{g} \cong ad(\mathfrak{g}) \subset \mathrm{End}(\mathfrak{g})$. Therefore $G$ is the identity component of $\mathrm{Aut}^{\mathbf{C}^*}(X, \omega)$ and $Lie (G) = \mathfrak{g}$.   
Q.E.D. 
\vspace{0.2cm}

{\bf Proposition 3} {\em The symplectic variety $X$ coincides with the closure of a 
coadjoint orbit of $\mathfrak{g}^*$.} 
\vspace{0.2cm}

{\em Proof}. 
Since $G$ is the identity component of $\mathrm{Aut}^{{\mathbf C}^*}(X, \omega)$, 
$X$ is stable under the coadjoint action of $G$ on $\mathfrak{g}$. Hence $X$ is a 
union of $G$-orbits. The $G$-orbits in $X$ are symplectic leaves of the Poisson variety $X$. 
In our case, since $X$ has only symplectic singularities, $X$ has only finitely many symplectic leaves by [Ka]. Therefore $X$ consists of finite number of $G$-orbits; hence there is an open dense $G$-orbit and $X$ is the closure of such an orbit. Q.E.D.    
\vspace{0.2cm}

For the unipotent radical $U$ of $G$, let us denote by $\mathfrak{n}$ its Lie algebra \footnote{The ideal $\mathfrak{n}$ is actually the nilradical of $\mathfrak{g}$ when 
$\mathfrak{g}$ has trivial center.}.   
Assume that $\mathfrak{n} \ne 0$. Then the center $z(\mathfrak{n})$ of $\mathfrak{n}$ is also non-trivial because $\mathfrak{n}$ is a nilpotent Lie algebra. Moreover $z(\mathfrak{n})$ is an ideal of $\mathfrak{g}$. In fact, it is enough to prove that, if $y \in \mathfrak{g}$ and $z \in z(\mathfrak{n})$, then $[x, [y, z]] = 0$ for any $x \in \mathfrak{n}$. Consider the Jacobi identity 
$$ [x, [y,z]] + [y, [z, x]] + [z, [x,y]] = 0. $$ 
First, since $z \in z(\mathfrak{n})$, one has $[z, x] = 0$. Next, since $\mathfrak{n}$ is an ideal of $\mathfrak{g}$, we have $[x,y] \in \mathfrak{n}$; hence $[z, [x,y]] = 0$. It then follows from the Jacobi identity that $[x, [y, z]] = 0$.
\vspace{0.2cm}

{\bf Proposition 4}. {\em Let $\mathfrak{g}$ be a complex Lie algebra with trivial center whose adjoint group $G$ is a linear algebraic group. Assume that  
$\mathfrak{n} \ne 0$. Let $O$ be a coadjoint orbit of $\mathfrak{g}^*$ with the following properties} 

(i) {\em $O$ is preserved by the scalar $\mathbf{C}^*$-action on $\mathfrak{g}^*$}; 

(ii) {\em $T_0\bar{O} = \mathfrak{g}^*$, where $T_0\bar{O}$ denotes the tangent space of 
the closure $\bar{O}$ of $O$ at the origin}.

{\em Then  $\bar{O} - O$ contains infinitely many coadjoint orbits; in particular $\bar{O}$ has infinitely many symplectic leaves.}  

{\bf Remark}. This proposition shows that $\bar{O}$ cannot have symplectic singularities. In fact, if $\bar{O}$ has symplectic singularities, it has only finitely many symplectic leaves by [Ka].

{\em Proof of Proposition 4}.  By a result of Mostow [Mos] (cf. [Ho], VIII, Theorem 3.5, 
Theorem 4.3), $G$ is a semi-direct product of a reductive subgroup $L$ and the unipotent radical $U$. Therefore we have a decomposition  $\mathfrak{g} = \mathfrak{l} \oplus \mathfrak{n}$. Take an element $\phi \in O$. Then $\phi$ is a linear function on $\mathfrak{g}$, which restricts to a 
{\em non-zero} function on $z(\mathfrak{n})$. In fact, if $\phi$ is zero on $z(\mathfrak{n})$, then $O \subset (\mathfrak{g}/z(\mathfrak{n}))^*$ and hence $\bar{O} \subset  (\mathfrak{g}/z(\mathfrak{n}))^*$, which contradicts the assumption (ii). 
We put $\bar{\phi} := \phi \vert_{z(\mathfrak n)} \ne 0$. 

Notice that the adjoint group $G$ is the subgroup of $\mathrm{GL}(\mathfrak{g})$ generated by all  elements of the form  
$\exp (ad \; v)$ with  $v \in \mathfrak{g}$. If $v \in z(\mathfrak{n})$, then 
$\exp (ad \; v) = id + ad \; v$ because $(ad \; v)^2 = 0$ for $v \in z(\mathfrak{n})$. 
Let $Z(U)$ be the identity component  of the center of the unipotent radical $U$. 
Then one can write $Z(U) = 1 + ad \; z(\mathfrak{n})$. 
By the assumption, the map $ad : \mathfrak{g} \to \mathrm{End}(\mathfrak{g})$ 
is an injection. Now we identify $z(\mathfrak{n})$ with $ad \; z(\mathfrak{n})$; then   
one can write $Z(U) = 1 + z(\mathfrak{n})$ and its group law is defined by $(1 + v)(1 + v') := 1 + (v + v')$ for $v, v' \in z(\mathfrak{n})$. Fix an element $v \in z(\mathfrak{n})$ and consider $Ad^*_{1 + v}(\phi)  \in 
\mathfrak{g}^*$. Since $(1+v)^{-1} = 1 -v$, the adjoint action $$Ad_{(1+v)^{-1}}:  \mathfrak{l} \oplus \mathfrak{n} \to 
\mathfrak{l} \oplus \mathfrak{n}$$ is defined by $$ x \oplus y \to x \oplus  (-[v, x] + y)$$ because 
$(ad \; v) (y) = 0$. By definition 
$$Ad^*_{1 + v}(\phi) (x \oplus y) = \phi (Ad_{(1+v)^{-1}}(x \oplus y)) = \phi (x \oplus y) - \bar{\phi}([v, x]).$$ Here notice that 
$[v, x] \in z(\mathfrak{n})$ and hence $\phi ([v,x]) =  \bar{\phi}([v, x])$. 

Since $z(\mathfrak{n})$ is an $\mathfrak{l}$-module by the Lie bracket, we decompose it into irreducible factors $z(\mathfrak{n}) = \bigoplus V_i$. Notice that it is the same as the irreducible decomposition of $z(\mathfrak{n})$ as a $[\mathfrak{l}, \mathfrak{l}]$-module if 
$[\mathfrak{l}, \mathfrak{l}] \ne 0$. In fact, the reductive Lie algebra $\mathfrak{l}$ 
is written as a direct sum of the semi-simple part and the center: $\mathfrak{l} = [\mathfrak{l}, \mathfrak{l}] \oplus z(\mathfrak{l})$. Since $z(\mathfrak{l})$ is an Abelian Lie algebra, 
$z(\mathfrak{n})$ can be written as a direct sum $\oplus V_{\alpha}$ of the weight spaces for 
$z(\mathfrak{l})$. 
The semisimple part $[\mathfrak{l}, \mathfrak{l}]$ acts on  each weight space $V_{\alpha}$; 
hence $V_{\alpha}$ is a direct sum of irreducible $[\mathfrak{l}, \mathfrak{l}]$ modules. 
These  irreducible $[\mathfrak{l}, \mathfrak{l}]$ modules are stable under the $z(\mathfrak{l})$-action and, hence are irreducible $\mathfrak{l}$-modules.

When $\dim V_i = 1$ for some $i$, this $V_i$ is an ideal of $\mathfrak{g}$. By the assumption (i), one can write $\bar{O} = \mathrm{Spec} R$ with a graded $\mathbf{C}$-algebra $R = 
\oplus_{j \ge 0} R_j$.  By (ii) we see that $R_1 = \mathfrak{g}$. 
Take a generator $x$ of a 1-dimensional space $V_i$. Then $x$ generates a Poisson ideal $I$ of $R$ and $Y := \mathrm{Spec}(R/I)$ is a closed Poisson subscheme of $\bar{O}$ of codimension $1$. Moreover, $Y$ is stable under the $G$-action. Since $\dim Y$ is odd, $Y$ contains infinitely many coadjoint orbits.   

In the remainder we assume that $\dim V_i > 1$ for all $i$.  In this case $[\mathfrak{l}, \mathfrak{l}] \ne 0$. 
Since $\bar{\phi} \ne 0$, we can choose an $i$ such that $\phi \vert_{V_i} \ne 0$. 
We fix a Cartan subalgebra $\mathfrak{h}$ of the semisimple Lie algebra $[\mathfrak{l}, \mathfrak{l}]$ and 
choose a set $\Delta$ of simple roots from the root system $\Phi$. We define 
$\mathfrak{n}^+ := \bigoplus_{\alpha \in \Phi^+}[\mathfrak{l}, \mathfrak{l}]_{\alpha}$.  
Let $v_0 \in V_i$ be a highest weight vector of the irreducible $[\mathfrak{l}, \mathfrak{l}]$-module $V_i$. Then one has $[v_0, \mathfrak{n}^+] = 0$ and, in particular, $\phi ([v_0, \mathfrak{n}^+] ) = 0$. Moreover, we may assume that 
$\bar{\phi} (v_0) \ne 0$ by replacing $\phi$ by a suitable $Ad^*_g(\phi)$ with $g \in L$. 
This is possible. In fact, if $Ad^*_{g}(\bar{\phi})(v_0) = 0$ for all $g$, then $\bar{\phi}$ is zero on the vector subspace of $V_i$ spanned by all $Ad_g (v_0)$. But, since $V_i$ is an irreducible  $L$-representation, such a subspace coincides with $V_i$. This contradicts the fact that $\bar{\phi}\vert_{V_i} \ne 0$.  
Since $v_0$ is a highest weight vector of a non-trivial $[\mathfrak{l}, \mathfrak{l}]$-irreducible 
module $V_i$, $[v_0, h]$ is a multiple of $v_0$ by a non-zero constant for an $h \in \mathfrak{h}$.  
Since $\bar{\phi} (v_0) \ne 0$, we also have $\bar{\phi}([v_0, h]) \ne 0$ for this $h \in \mathfrak{h}$. 
 
Let us consider $\bar{\phi}_{v_0} := \bar{\phi}([v_0, \cdot ]) \vert_{[\mathfrak{l}, \mathfrak{l}]}$.  By definition $\bar{\phi}_{v_0}$ 
is an element of $[\mathfrak{l}, \mathfrak{l}]^*$. By the Killing form it is identified with an 
element of $[\mathfrak{l}, \mathfrak{l}]$. The two facts $\bar{\phi}_{v_0}(\mathfrak{n}^+) = 0$ and 
$\bar{\phi}_{v_0}(h) \ne 0$ mean that $\bar{\phi}_{v_0}$ is {\em not} a nilpotent element 
of $[\mathfrak{l}, \mathfrak{l}]$.

For such $v_0$ and $\phi$, we consider 
$Ad^*_{1 + t^{-1}v_0}(t \phi)$, 
with $t \in \mathbf{C}^*$. One can write 
$$ Ad^*_{1 + t^{-1}v_0}(t \phi)(x \oplus y) = t\phi (x \oplus y) - \bar{\phi}([v_0, x]). $$ 
Thus one has  
$$  \lim_{t \to 0}Ad^*_{1 + t^{-1}v_0}(t \phi)(x \oplus y) = -\bar{\phi}([v_0, x]). $$
By definition  $Ad^*_{1 + t^{-1}v_0}(t \phi) \in O$. Thus $\lim_{t \to 0}Ad^*_{1 + t^{-1}v_0}(t \phi) \in \bar{O}$. Moreover, by the equality above, we see that  
$\lim_{t \to 0}Ad^*_{1 + t^{-1}v_0}(t \phi)\vert_{\mathfrak{n}} = 0$; thus it can be regarded as an element of $(\mathfrak{g}/\mathfrak{n})^* = \mathfrak{l}^*$.
 
Furthermore, we have $$\lim_{t \to 0}Ad^*_{1 + t^{-1}v_0}(t \phi)\vert_{[\mathfrak{l}, \mathfrak{l}]} = -\bar{\phi}_{v_0}, $$ which can be regarded as an element of $[\mathfrak{l}, \mathfrak{l}]$ by the identification $[\mathfrak{l}, \mathfrak{l}]^* \cong [\mathfrak{l}, \mathfrak{l}]$. As remarked above, this is not a nilpotent element.  
 
Let us write $\mathfrak{l}$ as a direct sum of the semi-simple part and the center: $\mathfrak{l} = [\mathfrak{l}, \mathfrak{l}] \oplus z(\mathfrak{l})$. There is an $L$-equivariant isomorphism $\mathfrak{l}^* \cong  [\mathfrak{l}, \mathfrak{l}]^* \oplus z(\mathfrak{l})^*$. Here $L$ acts trivially on the second factor $z(\mathfrak{l})^*$. Therefore, every coadjoint orbit of $\mathfrak{l}^*$ is a pair of a coadjoint orbit of $[\mathfrak{l}, \mathfrak{l}]^*$ and an element of $z(\mathfrak{l})^*$. 
 
In our situation, we can write $$\bar{\phi}([v_0, \cdot ]) = \bar{\phi}_{v_0} \oplus \bar{\phi}([v_0, \cdot ])\vert_{z(\mathfrak{l})}.$$ 
 
We can apply the same argument for $\lambda \phi$ with an arbitrary $\lambda \in 
\mathbf{C}^*$ to conclude that $\lambda \bar{\phi}([v_0, \cdot ]) \in \bar{O}$. 
One can write 
$$\lambda \cdot \bar{\phi}([v_0, \cdot ) = \lambda \bar{\phi}_{v_0} \oplus \lambda \cdot \bar{\phi}([v_0, \cdot ])\vert_{z(\mathfrak{l})}.$$

Since $\bar{\phi}_{v_0}$ is not an nilpotent element, we see that $\lambda \bar{\phi}_{v_0}$  $(\lambda \in \mathbf{C}^*)$ are contained in mutually different coadjoint orbits of $[\mathfrak{l}, \mathfrak{l}]^*$. 

Therefore, $\lambda \cdot \bar{\phi}([v_0, \cdot ])$ $(\lambda \in \mathbf{C}^*)$ are also contained in mutually different coadjoint orbits of 
$\mathfrak{l}^*$. Q.E.D. 

\vspace{0.2cm}

{\em Proof of Theorem}. We already know that $wt(\omega) = 1$ or $wt(\omega) = 2$. In the latter case $(X, \omega)$ is isomorphic to $(\mathbf{C}^{2d}, \omega_{st})$. So we assume that $wt(\omega) = 1$. By Propositions 2, 3 and 4, $(X, \omega)$ is isomorphic to 
a coadjoint orbit closure $(\bar{O}, \omega_{KK})$ of a complex reductive Lie algebra $\mathfrak{g}$ together with the Kirillov-Kostant form. Assume that $\mathfrak{g}$ is not semisimple, i.e., $\mathfrak{g}$ has non-trivial center $z(\mathfrak{g})$. 
In our case $O$ is preserved by the scalar $\mathbf{C}^*$-action on $\mathfrak{g}^*$. Such an orbit is contained in $(\mathfrak{g}/z(\mathfrak{g}))^*$. This contradicts the fact that 
$T_0\bar{O} = \mathfrak{g}^*$. Hence $\mathfrak{g}$ is semisimple. For a semisimple Lie algebra, a coadjoint orbit is identified with an adjoint orbit by the Killing form. A coadjoint orbit 
preserved by the scalar $\mathbf{C}^*$-action corresponds to a nilpotent orbit by this identification.  
\vspace{0.2cm}

\begin{center}
{\bf References}
\end{center}

[Be] Beauville,  A. : Symplectic singularities,  
Invent.  Math.  {\bf 139} (2000),  541-549

[Hi]{ Hinich,  V. : On the singularities of 
nilpotent orbits,  Israel J.  Math.  {\bf 73} (1991),  297-308

[Ho] Hochschild, G. P.: Basic theory of algebraic groups and Lie algebras, Graduate 
Texts in Mathematics {\bf 75} (1981)  

[Ka] Kaledin, D.: Symplectic varieties from the Poisson point of view, 
J. Reine Angew. Math. {\bf 600} (2006), 135-160

[K-P] Kraft, H., Procesi, C.: On the geometry of conjugacy classes in classical groups, Comment Math. Helv. {\bf 57} (1982), 539-602 

[Mos] Mostow, G.D.: Fully reducible subgroups of algebraic groups, Amer. J. Math. {\bf 78} (1956) 200-221

[Na 1] Namikawa, Y.: A finiteness theorem on symplectic singularities, Compositio Math. {\bf 152} (2016), 1225-1236  

[Na 2] Namikawa, Y.: On the structure of homogeneous symplectic varieties of complete intersection, Invent. Math. {\bf 193} (2013) 159-185 

[Na 3] Namikawa, Y.: Equivalence of symplectic singularities, Kyoto Journal of Mathematics, 
{\bf 53}, No.2 (2013), 483-514 

[Na 4] Namikawa, Y.: Flops and Poisson deformations of symplectic varieties, Publ. Res. Inst. Math. Sci. {\bf 44} (2008), 259 - 314  

[Pa] Panyushev,  D. I. : Rationality of singularities and the 
Gorenstein property of nilpotent orbits,  Funct.  Anal.  Appl.  {\bf 25} 
(1991),  225-226

[Pro] Procesi, C.: Lie groups, an approach through invariants and representations, 
(2007), UTX, Springer 

\vspace{0.3cm}

Department of Mathematics, Graduate school of Science, Kyoto University 

e-mail: namikawa@math.kyoto-u.ac.jp 

\end{document}